\documentclass{amsart}

\usepackage[dvipsnames]{xcolor}
\usepackage{amsmath}
\usepackage{enumerate}
\usepackage{enumitem}
\usepackage{amssymb}
\usepackage{extarrows}
\usepackage{amsthm}
\usepackage{thmtools}
\usepackage{mathrsfs}
\usepackage{tikz-cd}
\usepackage[colorlinks, linkcolor = SeaGreen, citecolor = SeaGreen, urlcolor = SeaGreen]{hyperref}

\hypersetup{colorlinks = true}

\newtheorem{thm}{Theorem}[section]

\newtheorem{coro}[thm]{Corollary}

\theoremstyle{definition}

\theoremstyle{remark}

\title{Commuting Toeplitz and Little Hankel Operators on the Bidisk}

\author[Lei Chen]{Lei Chen}
\address{\textsuperscript{1} College of Mathematics and Statistics, Chongqing University, Chongqing, 401331, P. R. China}
\email{9733618@qq.com}

\keywords{Toeplitz operator; Pluriharmonic Hardy space; Bidisk; Little Hankel operator}

\subjclass[2020]{47A16, 47B35, 30H20}

\begin{document}
	\begin{abstract}
        Necessary and sufficient conditions are obtained for two Toeplitz operators and for a Toeplitz operator and a little Hankel operator to commute on the pluriharmonic Hardy space of the bidisk.
	\end{abstract}
	
	\maketitle
	
	\section{Introduction}\label{S1}

    Let $\mathbb{D}$ denote the open unit disk in the complex plane $\mathbb{C}$.
    Its boundary is the unit circle $\mathbb{T}$. The bidisk $\mathbb{D}^2$ and the torus $\mathbb{T}^2$
    are this subsets of $\mathbb{C}^2$ which are Cartesian products of two copies $\mathbb{D}$
    and $\mathbb{T}$, respectively. Let $\mathrm{d}\sigma(z)$ be the normalized
    Haar measure on $\mathbb{T}^2$. The Hardy space $H^2(\mathbb{D}^2)$ is the 
    closure of the analytic polynomials in $L^2(\mathbb{T}^2,\mathrm{d}\sigma)$
    (or $L^2(\mathbb{T}^2)$). For $u,v\in\mathbb{Z}$ let $\chi_{u,v}$ be the function on $\mathbb{T}^2$ defined by
    $$
    \chi_{u,v}(e^{i\theta_1},e^{i\theta_2})=e^{i(u\theta_1+v\theta_2)}.
    $$
    For $p=1,2,\infty$, we define the Hardy space:
    $$
    H^p(\mathbb{D})=\left\{f\in L^p(\mathbb{T}):\int_0^{2\pi}f(e^{i\theta})\chi_n(e^{i\theta})\mathrm{d}\theta=0\text{ for }n>0\right\}.
    $$
    Similarly, the Hardy space on the bidisk is defined by
    $$
    H^p(\mathbb{D}^2)=\left\{f\in L^p(\mathbb{T}^2):\int_0^{2\pi}\int_0^{2\pi} f(e^{i\theta_1},e^{i\theta_2})\chi_{u,v}(e^{i\theta_1},e^{i\theta_2})\mathrm{d}\theta_1\mathrm{d}\theta_2=0\text{ for }u>0\text{ or }v>0\right\},
    $$
    and $H^2(\mathbb{D}^2)$ is the closure in the $L^2$-norm of the analytic trigonometric polynomials $\mathscr{P}_+$.

    We define the pluriharmonic Hardy space $h^2(\mathbb{T}^2)$ by
    $$
    h^2(\mathbb{T}^2)=H^2(\mathbb{T}^2)+\overline{H^2(\mathbb{T}^2)}.
    $$
    $h^2(\mathbb{T}^2)$ is a Hilbert with the inner product inherited from $L^2(\mathbb{T}^2)$.
    Then it is easy to vertify that $\{\chi_{u,v}\}_{(u,v)\in S}$
    is an orthonormal basis of $h^2(\mathbb{T}^2)$, where 
    $$
    S=\{(u,v)\in\mathbb{Z}^2:u\geqslant0,v\geqslant0\}\cup\{(u,v)\in\mathbb{Z}^2:u\leqslant0,v\leqslant0\}.
    $$
    For $f\in L^2(\mathbb{T}^2)$, the Toeplitz operator $\widehat{T}_f$ with symbol $f$
    is densely defined on the pluriharmonic Hardy space $h^2(\mathbb{T}^2)$ by
    $$
    \widehat{T}_fh=Q(fh)
    $$
    for all polynomials $h$, where $Q$ is the orthogonal projection from $L^2(\mathbb{T}^2)$ onto $h^2(\mathbb{T}^2)$.
    Let $P$ denote the orthogonal projection from $L^2(\mathbb{T}^2)$ onto $H^2(\mathbb{T}^2)$.
    Then for each function $f$ in $L^2(\mathbb{T}^2)$, note
    $$
    \begin{aligned}
        &f_{++}=Pf,\\
        &f_{--}=(Q-P)f.
    \end{aligned}
    $$

    For $\varphi\in L^\infty(\mathbb{T}^2)$, the Hankel operator $\widehat{H}_\varphi$ with symbol $\varphi$ is defined on the pluriharmonic Hardy
    space $h^2(\mathbb{T}^2)$ by
    $$
    \widehat{H}_\varphi:h^2(\mathbb{T}^2)\to\left(h^2(\mathbb{T}^2)\right)^\perp,\quad \widehat{H }_\varphi f=(I-Q)(\varphi f),\quad\forall f\in h^2(\mathbb{T}^2).
    $$

    For $f\in L^2(\mathbb{T}^2)$, define $\widetilde{f}(z_1,z_2)=f(\overline{z_1},\overline{z_2})$. Define the operator $U$ on $L^2(\mathbb{T}^2)$ by
    $$
    Uf(z_1,z_2)=\widetilde{f}(z_1,z_2).
    $$
    The little Hankel operator $\widehat{\Gamma}_\psi$ with symbol $\psi\in L^\infty(\mathbb{T}^2)$ on $h^2(\mathbb{T}^2)$ is defined by
    $$
    \widehat{\Gamma}_\psi f=QU(\psi f).
    $$

    The commutativity of Toeplitz operators is a classical problem
    in operator theory. On the Hardy space of the unit disk,
    Brown and Halmos~\cite{BH1964} proved in 1964 that two Toeplitz
    operators with bounded symbols commute if and only if both
    symbols are analytic, both are co-analytic, or a nontrivial
    linear combination of the symbols is constant. They also
    studied when a product of Toeplitz operators is itself a
    Toeplitz operator. These questions, commonly referred to as
    Brown--Halmos type problems, have subsequently been studied
    on other function spaces. On the Bergman space of the unit
    disk, Axler and \v{C}u\v{c}kovi\'{c}~\cite{AC1991}
    characterized commuting Toeplitz operators with harmonic
    symbols in 1991, using a conformally invariant mean-value
    characterization of harmonic functions. Ahern and
    \v{C}u\v{c}kovi\'{c}~\cite{AC2001} subsequently obtained
    Brown--Halmos type results for products with harmonic
    symbols in 2001, under regularity assumptions on the
    resulting symbol, including boundedness of its invariant
    Laplacian. They also settled the zero-product problem for
    harmonic symbols and constructed counterexamples for more
    general symbols. These developments show that the symbol
    class and the underlying projection both play a role in
    extending the classical algebraic results.

    Toeplitz product problems on analytic Hardy spaces in several
    variables have been investigated through zero products,
    semi-commutators, and finite sums of products. On the Hardy
    space of the bidisk, Gu and Zheng~\cite{GZ1997} characterized
    when a semi-commutator vanishes in 1997. They also proved
    that a finite-rank semi-commutator must be zero and gave
    examples of nonzero compact semi-commutators with continuous
    symbols. On the Hardy space of the polydisk,
    Ding~\cite{Ding2003} studied zero products, finite-rank
    products, and products that remain Toeplitz in 2003.
    In particular, a product of two Toeplitz operators has
    finite rank precisely when it is zero, which occurs
    precisely when one of the factors is zero. Choe, Koo and
    Lee~\cite{CKL2011} considered finite sums of products with
    pluriharmonic symbols on the Hardy space over the unit ball
    in 2011 and investigated their compactness and finite rank.
    These works address both exact product identities and the
    extent to which such identities can hold modulo compact
    or finite-rank operators.

    For commuting Toeplitz operators on analytic Hardy spaces,
    Ding, Sun and Zheng~\cite{DSZ2012} obtained necessary and
    sufficient conditions for bounded symbols on the bidisk in
    2012. Their main tools were the Berezin transform and
    harmonic extensions. Kong and Lu~\cite{KL2015} studied the
    corresponding problem on the polydisk in 2015 and showed
    that commutativity is equivalent to the separate harmonicity
    of the Berezin transform of the commutator. Both approaches
    relate operator equations to properties of functions on
    the underlying domain. The study of Hardy-Toeplitz operators
    on the bidisk also includes algebraic and spectral questions
    beyond commutativity. In 2023, Ding, Wang and
    Li~\cite{DWL2023} examined invariant subspaces, eigenvalues,
    eigenvectors and numerical ranges, providing further
    information about the relation between Toeplitz symbols
    and operator properties in this setting.

    On the pluriharmonic Hardy space $h^2(\mathbb{T}^2)$,
    the presence of both analytic and co-analytic Fourier modes
    leads to different commutation phenomena. Liu and
    Ding~\cite{LD2013} characterized commutativity and
    semi-commutativity for bounded analytic symbols in 2013.
    Their results show that two analytic symbols do not
    automatically give commuting Toeplitz operators on this
    space, since the action on the co-analytic part must also
    be considered. Ding and Sang~\cite{DS2017} obtained
    characterizations for bounded pluriharmonic symbols in
    2017. Their criteria involve the coordinate dependence
    of the analytic and co-analytic components and the
    existence of a constant linear combination of the symbols.
    In particular, the classification for symbols each
    depending on a single coordinate follows from their
    general commutativity theorem. Sang and Ding~\cite{SD2018}
    subsequently characterized when an analytic Toeplitz
    operator commutes with a co-analytic Toeplitz operator
    on the same space in 2018. Zhang and Dong~\cite{ZD2020}
    studied monomial symbols in 2020 and completely
    characterized finite-rank commutators and semi-commutators
    on the pluriharmonic Hardy spaces of the torus and the
    unit sphere. They also obtained nontrivial commuting
    and semi-commuting examples. These results make the
    coordinate dependence and Fourier support of the symbols
    central to the study of commutativity on pluriharmonic
    Hardy spaces.

    The mixed commutation problem for Toeplitz and Hankel
    operators forms a related line of research.
    Mart\'inez-Avenda\~{n}o~\cite{MA2000} completely classified
    commuting Toeplitz and Hankel operators on the classical
    Hardy space in 2000. His description of nontrivial pairs
    involves relations between the symbols and certain subsets
    of the unit circle. In 2026, Zou and Yu~\cite{ZY2026}
    obtained necessary and sufficient conditions for a Toeplitz
    operator with a bounded pluriharmonic symbol to commute
    with a Hankel operator with a bounded symbol on the Hardy
    space of the bidisk. These results concern operators on
    analytic Hardy spaces. They provide a point of comparison
    for the mixed problem on $h^2(\mathbb{T}^2)$, where the
    projection retains both analytic and co-analytic Fourier
    modes. In this setting, it is natural to investigate how
    the Fourier coefficients of the symbols determine
    commutation with the little Hankel operator defined above.

    In this paper, we study Toeplitz--Toeplitz and
    Toeplitz--little Hankel commutation relations on
    $h^2(\mathbb{T}^2)$ by a direct Fourier coefficient method.
    For symbols each depending on a single coordinate,
    Theorem~2.1 gives a Fourier coefficient proof of the
    Toeplitz classification discussed above. Theorem~3.1
    characterizes commutation with the little Hankel operator
    $\widehat{\Gamma}_\psi$ defined above. Apart from the cases
    $\psi=0$ or $\varphi$ constant, commutativity holds if and
    only if $\tilde{\varphi}=\varphi$ and
    $\psi=c\varphi+d$ for some constants $c,d\in\mathbb{C}$.
    For arbitrary $\varphi,\psi\in L^\infty(\mathbb{T}^2)$,
    Theorems~2.2 and~3.2 give necessary and sufficient Fourier
    coefficient criteria using $S$-restricted convolutions.
    These criteria require neither analyticity nor
    pluriharmonicity of the symbols and allow infinitely many
    nonzero Fourier coefficients. The restrictions on the
    summation indices explicitly encode the intermediate
    projection $Q$. The coordinate reflection $U$ gives
    \[
    \widehat{T}_\varphi\widehat{\Gamma}_\psi
    =
    \widehat{\Gamma}_\psi\widehat{T}_\varphi
    \quad\Longleftrightarrow\quad
    \widehat{T}_{\tilde{\varphi}}\widehat{T}_\psi
    =
    \widehat{T}_\psi\widehat{T}_\varphi.
    \]
    Consequently, both commutation problems can be treated
    using the same Toeplitz product coefficient formulas.
    We also establish in Corollary~2.3 the product identity
    \[
    \widehat{T}_\varphi\widehat{T}_\psi
    =
    \widehat{T}_{\varphi\psi}
    =
    \widehat{T}_\psi\widehat{T}_\varphi
    \]
    when both Fourier supports are contained in the closed
    fourth quadrant. This sufficient condition permits
    mixed-sign frequencies and does not impose a finite-support
    assumption. The resulting formulation complements the
    structural classifications discussed above by giving
    explicit coefficient criteria for general bounded symbols
    and a common method for the two commutation problems.

    \section{Commutativity of Toeplitz Operators}

    In this section, we characterize the commutativity of two Toeplitz operators on the pluriharmonic Hardy space of the bidisc. We begin with the case of bounded symbols depending on a single variable.

    \begin{thm}\label{thm3.10}
        Let $\varphi,\psi\in L^\infty(\mathbb{T})$. $\widehat{T}_\varphi$ commutes with $\widehat{T}_\psi$ on $h^2(\mathbb{T}^2)$ if and only if one of the following conditions holds:

        \begin{enumerate}[label = \textup{(\alph*)}]
			\item one of $\varphi$ and $\psi$ is constant.
            \item $\varphi$ and $\psi$ are in the same variable $z_1$ or $z_2$:
            \begin{enumerate}[label = \textup{($\mathrm{b}_{\arabic*}$)}]
                \item $\varphi$ and $\psi$ are in $H^\infty(\mathbb{T})$; or
                \item $\varphi$ and $\psi$ are in $\overline{H^\infty(\mathbb{T})}$; or
                \item $\varphi=c\psi+d$, where $c$ and $d$ are constants.
            \end{enumerate}
            \item $\varphi$ and $\psi$ are in different variable $z_1$ and $z_2$: one of $\varphi$ and $\psi$ is in $H^\infty(\mathbb{T})$, and the other is in $\overline{H^\infty(\mathbb{T})}$.
		\end{enumerate}
    \end{thm}

    \begin{proof}[Proof]
        Since $\{\chi_{k,j}\}_{(k,j)\in S}$
        is an orthonormal basis of $h^2(\mathbb{T}^2)$, we only need consider
        $$
        \widehat{T}_\varphi\widehat{T}_\psi\chi_{k,j}=\widehat{T}_\psi\widehat{T}_\varphi\chi_{k,j},
        $$
        for all $(k,j)\in S$.

        If one of $\varphi$ and $\psi$ is constant, they commute obviously. Then we consider $\varphi$ and $\psi$ are both nonconstant.

        If $\varphi$ and $\psi$ are one variable function in $L^2(\mathbb{T}^2)$, firstly we assume that $\varphi$ and $\psi$ are variable in $z_1$,
        i.e.
        $$
        \varphi(z)=\sum_{n=-\infty}^{+\infty}\widehat{\varphi}(n)z_1^n,\quad\psi(z)=\sum_{m=-\infty}^{+\infty}\widehat{\psi}(m)z_1^m,
        $$
        where $\widehat{\varphi}(n)$ and $\widehat{\psi}(m)$ are the Fourier coefficients of $\varphi$ and $\psi$ in the $L^2(\mathbb T^2)$ sense respectively.
        Since $\widehat{T}_\varphi\widehat{T}_\psi\chi_{k,j}=\widehat{T }_\psi\widehat{T}_\varphi\chi_{k,j}$ for all $(k,j)\in S$, 
        we consider $k\geqslant$ and $j\geqslant0$ firstly. By calculating, we can get
        $$
        \widehat{T}_\varphi\widehat{T}_\psi\chi_{k,j}=Q[\varphi Q(\psi\chi_{k,j})]=\begin{cases}
            \displaystyle\sum_{m+k\geqslant0,i\geqslant0}\widehat{\varphi}(i-m-k)\widehat{\psi}(m)z_1^iz_2^j,\quad&k\geqslant0,j>0,\\
            \vspace{2mm}\\
            \displaystyle\sum_{n=-\infty}^{+\infty}\sum_{m=-\infty}^{+\infty}\widehat{\varphi}(n)\widehat{\psi}(m)z_1^{m+n+k},&k>0,j=0
        \end{cases}
        $$
        and
        $$
        \widehat{T}_\psi\widehat{T}_\varphi\chi_{k,j}=Q[\psi Q(\varphi\chi_{k,j})]=\begin{cases}
            \displaystyle\sum_{n+k\geqslant0,i\geqslant0}\widehat{\psi}(i-n-k)\widehat{\varphi}(n)z_1^iz_2^j,\quad&k\geqslant0,j>0,\\
            \vspace{2mm}\\
            \displaystyle\sum_{m=-\infty}^{+\infty}\sum_{n=-\infty}^{+\infty}\widehat{\varphi}(n)\widehat{\psi}(m)z_1^{m+n+k},&k>0,j=0.
        \end{cases}
        $$
        Obviously when $k>0,j=0$, $\varphi$ and $\psi$ are arbitrary in $L^2(\mathbb{T}^2)$. Then for any $k\geqslant0,j>0$, we get
        $$
        \sum_{m+k\geqslant0,i\geqslant0}\widehat{\varphi}(i-m-k)\widehat{\psi}(m)z_1^iz_2^j=\sum_{n+k\geqslant0,i\geqslant0}\widehat{\psi}(i-n-k)\widehat{\varphi}(n)z_1^iz_2^j,
        $$
        then fixing any $i\geqslant0$, we have
        $$
        \sum_{m+k\geqslant0}\widehat{\varphi}(i-m-k)\widehat{\psi}(m)=\sum_{n+k\geqslant0}\widehat{\psi}(i-n-k)\widehat{\varphi}(n).
        $$
        Let $l=m+k$ in left-hand side and $l=n+k$ in right-hand side,
        \begin{equation}\label{3.10.1}
            \sum_{l\geqslant0}\widehat{\varphi}(i-l)\widehat{\psi}(l-k)=\sum_{l\geqslant0}\widehat{\psi}(i-l)\widehat{\varphi}(l-k).
        \end{equation}
        Take $i+1$ and $k+1$ in (\ref{3.10.1}),
        \begin{equation}\label{3.10.2}
            \sum_{l\geqslant0}\widehat{\varphi}(i+1-l)\widehat{\psi}(l-k-1)=\sum_{l\geqslant0}\widehat{\psi}(i+1-l)\widehat{\varphi}(l-k-1).
        \end{equation}
        Replace $l$ by $l-1$ in (\ref{3.10.1}),
        \begin{equation}\label{3.10.3}
            \sum_{l\geqslant1}\widehat{\varphi}(i-l+1)\widehat{\psi}(l-1-k)=\sum_{l\geqslant1}\widehat{\psi}(i-l+1)\widehat{\varphi}(l-1-k).
        \end{equation}
        Let (\ref{3.10.2}) minus (\ref{3.10.3}),
        $$
        \widehat{\varphi}(i+1)\widehat{\psi}(-(k+1))=\widehat{\psi}(i+1)\widehat{\varphi}(-(k+1)).
        $$
        Since $i\geqslant0$ and $k\geqslant0$, we have
        \begin{equation}\label{3.10.4}
            \widehat{\varphi}(p)\widehat{\psi}(q)=\widehat{\psi}(p)\widehat{\varphi}(q),\quad\forall p\geqslant1,q\leqslant-1.
        \end{equation}
        
        Similarly as $k\leqslant$ and $j\leqslant0$, we have
        \begin{equation}\label{3.10.5}
            \widehat{\varphi}(p)\widehat{\psi}(q)=\widehat{\psi}(p)\widehat{\varphi}(q),\quad\forall p\leqslant-1,q\geqslant1.
        \end{equation}
        By (\ref{3.10.4}) and (\ref{3.10.5}), we can get
        \begin{equation}\label{3.10.6}
            \widehat{\varphi}(p)\widehat{\psi}(q)=\widehat{\psi}(p)\widehat{\varphi}(q),\quad\forall p,q\neq0,pq<0.
        \end{equation}

        Now we define $\psi_+=P\psi-\widehat{\psi}(0)$ and $\psi_-=\psi_{--}$, and similarly as $\varphi$.
        If $\psi_+\neq0$ and $\psi_-\neq0$, we can find $q_1<0$ s.t. $\widehat{\psi}(q_1)\neq0$, then for all $p>0$,
        $$
        \widehat{\varphi}(p)=\frac{\widehat{\varphi}(q_1)}{\widehat{\psi}(q_1)}\widehat{\psi}(p).
        $$
        Let 
        $$
        c=\frac{\widehat{\varphi}(q_1)}{\widehat{\psi}(q_1)}.
        $$
        We also can find $q_2>0$ s.t. $\widehat{\psi}(q_2)\neq0$. Then by (\ref{3.10.6}), for any $p<0$ we have
        $$
        \widehat{\varphi}(p)\widehat{\psi}(q_2)=\widehat{\psi}(p)\widehat{\varphi}(q_2).
        $$
        Hence
        $$
        \widehat{\varphi}(p)=\frac{\widehat{\varphi}(q_2)}{\widehat{\psi}(q_2)}\widehat{\psi}(p).
        $$
        Let
        $$
        c'=\frac{\widehat{\varphi}(q_2)}{\widehat{\psi}(q_2)}.
        $$
        Applying (\ref{3.10.6}) with $p=q_1<0$ and $q=q_2>0$, we have
        $$
        \widehat{\varphi}(q_1)=\frac{\widehat{\varphi}(q_2)}{\widehat{\psi}(q_2)}\widehat{\psi}(q_1),
        $$
        i.e.
        $$
        c=\frac{\widehat{\varphi}(q_1)}{\widehat{\psi}(q_1)}=\frac{\widehat{\varphi}(q_2)}{\widehat{\psi}(q_2)}=c'.
        $$
        Hence we can get $\widehat{\varphi}(n)=c\widehat{\psi}(n)$ for $n\neq0$, i.e. $\varphi=c\psi+d$, where $c$ and $d$ are constants. ($b_3$) holds.
        
        If $\psi_-=0$ and $\psi_+\neq0$ i.e. $\psi$ is in $H^\infty(\mathbb{T})$, by (\ref{3.10.6}), let $p<0,q>0$, we have
        $$
        \widehat{\varphi}(p)\widehat{\psi}(q)=0.
        $$
        Since $\psi_+\neq0$, there exists $q>0$ s.t. $\widehat{\psi}(q)\neq0$. It follows that $\widehat{\varphi}(p)=0$ for all $p<0$.
        Hence $\varphi_-=0$, i.e. $\varphi$ is in $H^\infty(\mathbb{T})$. ($b_1$) holds.

        If $\psi_-\neq0$ and $\psi_+=0$, similarly argument as above case and hence we have $\varphi$ and $\psi$ are in $\overline{H^\infty(\mathbb{T})}$. ($b_2$) holds.

        For $\varphi$ and $\psi$ depending on $z_2$, the argument proceeds in exactly the same way as above. Therefore $(b)$ holds.
        
        If $\varphi$ and $\psi$ depend on different variable, we may assume that $\varphi$ depends on $z_1$ and $\psi$ depends on $z_2$, i.e.
        $$
        \varphi(z)=\sum_{n=-\infty}^{+\infty}\widehat{\varphi}(n)z_1^n,\quad\psi(z)=\sum_{m=-\infty}^{+\infty}\widehat{\psi}(m)z_2^m.
        $$
        For $k>0,j=0$, we have
        $$
        \widehat{T}_\varphi\widehat{T}_\psi\chi_{k,0}=\sum_{n+k\geqslant0}\sum_{m\geqslant0}\widehat{\varphi}(n)\widehat{\psi}(m)z_1^{k+n}z_2^m+\sum_{n+k<0}\widehat{\varphi}(n)\widehat{\psi}(0)z_1^{k+n}
        $$
        and
        $$
        \widehat{T}_\psi\widehat{T}_\varphi\chi_{k,0}=\sum_{k+n\geqslant0}\sum_{m\geqslant0}\widehat{\varphi}(n)\widehat{\psi}(m)z_1^{k+n}z_2^m+\sum_{k+n\leqslant0}\sum_{m\leqslant0}\widehat{\varphi}(n)\widehat{\psi}(m)z_1^{k+n}z_2^m-\widehat{\varphi}(-k)\widehat{\psi}(0).
        $$
        By $\widehat{T}_\varphi\widehat{T}_\psi\chi_{k,0}=\widehat{T}_\psi\widehat{T}_\varphi\chi_{k,0}$, we have
        $$
        \sum_{k+n\leqslant0,m<0}\widehat{\varphi}(n)\widehat{\psi}(m)z_1^{k+n}z_2^m=0.
        $$
        Since $k+n\leqslant0$ and $k>0$, we have $n<0$. It follows that
        \begin{equation}\label{3.10.7}
            \widehat{\varphi}(n)\widehat{\psi}(m)=0,\quad\forall n<0,m<0.
        \end{equation}
        Next, using the same method, we obtain: for $k<0,j=0$,
        \begin{equation}\label{3.10.8}
            \widehat{\varphi}(n)\widehat{\psi}(m)=0,\quad\forall n>0,m>0;
        \end{equation}
        for $k=0,j>0$,
        $$
        \widehat{\varphi}(n)\widehat{\psi}(m)=0,\quad\forall n<0,m<0;
        $$
        for $k=0,j<0$,
        $$
        \widehat{\varphi}(n)\widehat{\psi}(m)=0,\quad\forall n>0,m>0.
        $$
        For $k>0,j>0$,
        $$
        \widehat{T}_\varphi\widehat{T}_\psi\chi_{k,j}=\sum_{k+n\geqslant0}\sum_{j+m\geqslant0}\widehat{\varphi}(n)\widehat{\psi}(m)z_1^{k+n}z_2^{j+m}+\widehat{\psi}(-j)\sum_{k+n<0}\widehat{\varphi}(n)z_1^{k+n}
        $$
        and
        $$
        \widehat{T}_\psi\widehat{T}_\varphi\chi_{k,j}=\sum_{k+n\geqslant0}\sum_{j+m\geqslant0}\widehat{\varphi}(n)\widehat{\psi}(m)z_1^{k+n}z_2^{j+m}+\widehat{\varphi}(-k)\sum_{j+m<0}\widehat{\psi}(m)z_2^{j+m}.
        $$
        Since $-j<0,n<-k<0$ and $-k<0,m<-j<0$, by (\ref{3.10.7}), both additional terms vanish. Hence the two expressions coincide.
        And for $k<0,j<0$, the argument is similar. Lastly, for $k=j=0$,
        $$
        \widehat{T}_\varphi\widehat{T}_\psi\chi_{0,0}=Q(\varphi\psi)=Q(\psi\varphi)=\widehat{T}_\psi\widehat{T}_\varphi\chi_{0,0}.
        $$
        Therefore, since $\varphi$ and $\psi$ are both nonconstant, by (\ref{3.10.7}), (\ref{3.10.8}),
        we have one of $\varphi$ and $\psi$ is in $H^\infty(\mathbb{T})$, the other is in $\overline{H^\infty(\mathbb{T})}$.

        Conversely, for condition ($a$), obviously.

        For condition $(b_1)$, we need to show
        $$
        \widehat{T}_\varphi\widehat{T}_\psi\chi_{k,j}=\widehat{T}_\psi\widehat{T}_\varphi\chi_{k,j},
        $$
        for all $(k,j)\in S$. Since $\varphi$ and $\psi$ are in $H^\infty(\mathbb{T})$, it is obviously for $k\geqslant0$ and $j\geqslant0$.
        Then we only need to consider $k\leqslant0$ and $j\leqslant0$. We may assume that $\varphi$ and $\psi$ are variable in $z_1$. For $k\leqslant0,j<0$,
        $$
        \widehat{T}_\varphi\widehat{T}_\psi\chi_{k,j}=\sum_{n+m+k\leqslant0,m\geqslant0,n\geqslant0}\widehat{\varphi}(n)\widehat{\psi}(m)z_1^{m+n+k}z_2^j=\widehat{T}_\psi\widehat{T}_\varphi\chi_{k,j},
        $$
        and for $k<0,j=0$,
        $$
        \widehat{T}_\varphi\widehat{T}_\psi\chi_{k,0}=\sum_{m\geqslant0,n\geqslant0}\widehat{\varphi}(n)\widehat{\psi}(m)z_1^{n+m+k}=\widehat{T}_\psi\widehat{T}_\varphi\chi_{k,0}.
        $$
        And similarly as the case that $\varphi$ and $\psi$ are variable in $z_2$. Hence $\widehat{T}_\varphi$ and $\widehat{T}_\psi$ commute.

        For condition $(b_2)$, using the same method as $(b_1)$.

        For condition $(b_3)$, we have
        $$
        \widehat{T}_\varphi\widehat{T }_\psi=(c\widehat{T }_\psi+dI)\widehat{T}_\psi=c\widehat{T}^2_\psi+d\widehat{T }_\psi=\widehat{T }_\psi(c\widehat{T }_\psi+dI)=\widehat{T }_\psi\widehat{ T }_\varphi.
        $$

        For condition $(c)$, assume that $\varphi$ is analytic in $z_1$ and $\psi$ is co-analytic in $z_2$, i.e. $\widehat{\varphi}(n)=0,n<0$ and $\widehat{\psi}(m)=0,m>0$.
        For $k>0,j>0$,
        $$
        \widehat{T}_\varphi\widehat{T}_\psi\chi_{k,j}=\sum_{n\geqslant0}\sum_{-j\leqslant m\leqslant0}\widehat{\varphi}(n)\widehat{\psi}(m)z_1^{k+n}z_2^{j+m}=\widehat{T}_\psi\widehat{T}_\varphi\chi_{k,j},
        $$
        and for $k<0,j<0$, 
        $$
        \widehat{T}_\varphi\widehat{T}_\psi\chi_{k,j}=\sum_{0\leqslant n\leqslant-k}\sum_{m\leqslant0}\widehat{\varphi}(n)\widehat{\psi}(m)z_1^{k+n}z_2^{j+m}=\widehat{T}_\psi\widehat{T}_\varphi\chi_{k,j}.
        $$
        Then we consider the boundary case. For $k>0,j=0$,
        $$
        \widehat{T}_\varphi\widehat{T}_\psi\chi_{k,0}=\widehat{\psi}(0)\sum_{n\geqslant0}\widehat{\varphi}(n)z_1^{n+k}=\widehat{T}_\psi\widehat{T}_\varphi\chi_{k,0};
        $$
        for $k<0,j=0$,
        $$
        \widehat{T}_\varphi\widehat{T}_\psi\chi_{k,0}=\sum_{m\leqslant0,n\geqslant0,n+k\leqslant0}\widehat{\varphi}(n)\widehat{\psi}(m)z_1^{k+n}z_2^m+\widehat{\psi}(0)\sum_{n\geqslant0,k+n>0}\widehat{\varphi}(n)z_1^{k+n}=\widehat{T}_\psi\widehat{T}_\varphi\chi_{k,0};
        $$
        for $k=0,j>0$,
        $$
        \widehat{T}_\varphi\widehat{T}_\psi\chi_{0,j}=\sum_{m\leqslant0,n>0,m+j\geqslant0}\widehat{\varphi}(n)\widehat{\psi}(m)z_1^nz_2^{m+j}+\widehat{\varphi}(0)\sum_{m\leqslant0}\widehat{\psi}(m)z_2^{m+j}=\widehat{T}_\psi\widehat{T}_\varphi\chi_{0,j};
        $$
        for $k=0,j<0$,
        $$
        \widehat{T}_\varphi\widehat{T}_\psi\chi_{0,j}=\widehat{\varphi}(0)\sum_{m\leqslant0}\widehat{\psi}(m)z_2^{m+j}=\widehat{T}_\psi\widehat{T}_\varphi\chi_{0,j};
        $$
        for $k=j=0$,
        $$
        \widehat{T}_\varphi\widehat{T}_\psi\chi_{0,0}=Q(\varphi\psi)=Q(\psi\varphi)=\widehat{T}_\psi\widehat{T}_\varphi\chi_{0,0}.
        $$
        The same holds in all other cases for $\varphi$ and $\psi$. 
        
        Therefore, commutativity holds.
        This completes the proof.
    \end{proof}

    Now we define two "$S$-restricted convolutions" by
    $$
    (\widehat{\varphi}\ast_S^{(i,j)}\widehat{\psi})(M,N)=\sum_{(k+i,l+j)\in S}\widehat{\varphi}(M-i-k,N-j-l)\widehat{\psi}(k,l),\quad\forall(i,j)\in S
    $$
    and
    $$
    (\widehat{\psi}\ast_S^{(i,j)}\widehat{\varphi})(M,N)=\sum_{(m+i,n+j)\in S}\widehat{\psi}(M-i-m,N-j-n)\widehat{\varphi}(m,n),\quad\forall(i,j)\in S,
    $$
    where $\widehat{\varphi}(m,n)$ and $\widehat{\psi}(k,l)$ are the Fourier coefficients of $\varphi$ and $\psi$ in $L^\infty(\mathbb{T}^2)$ respectively. 
    Fixing $(i,j),(M,N)\in S$, by the Cauchy-Schwarz inequality and Parseval's identity,
    $$
    \begin{aligned}
        &\sum_{(k+i,l+j)\in S}\left|\widehat{\varphi}(M-i-k,N-j-l)\widehat{\psi}(k,l)\right|\\
        &\leqslant\left(\sum_{(k+i,l+j)\in S}|\widehat{\varphi}(M-i-k,N-j-l)|^2\right)^{\frac{1}{2}}\left(\sum_{(k+i,l+j)\in S}|\widehat{\psi}(k,l)|^2\right)^{\frac{1}{2}}\\
        &\leqslant\left(\sum_{k,l\in\mathbb{Z}}|\widehat{\varphi}(M-i-k,N-j-l)|^2\right)^{\frac{1}{2}}\left(\sum_{k,l\in\mathbb{Z}}|\widehat{\psi}(k,l)|^2\right)^{\frac{1}{2}}\\
        &=\Vert\varphi\Vert_{L^2(\mathbb{T}^2)}\Vert\psi\Vert_{L^2(\mathbb{T}^2)}<\infty.
    \end{aligned}
    $$
    The last equality follows because $(k,l)\mapsto(M-i-k,N-j-l)$ is a bijection of $\mathbb{Z}^2$.
    Thus $(\widehat{\varphi}\ast_S^{(i,j)}\widehat{\psi})(M,N)$ is absolutely convergent.
    Interchanging $\varphi$ and $\psi$ shows that the other restricted convolution is also absolutely convergent.
    We then have the following theorem for bounded symbols in two variables.

    \begin{thm}\label{thm3.11}
        Let $\varphi,\psi\in L^\infty(\mathbb{T}^2)$. $\widehat{T}_\varphi$ commutes with $\widehat{T}_\psi$ on $h^2(\mathbb{T}^2)$
        if and only if for any $(i,j)\in S$,
        $$
        \widehat{\varphi}\ast_S^{(i,j)}\widehat{\psi}=\widehat{\psi}\ast_S^{(i,j)}\widehat{\varphi}\quad\text{on $S$.}
        $$
    \end{thm}

    \begin{proof}[Proof]
        Suppose first that $\widehat{T}_\varphi$ commutes with $\widehat{T}_\psi$.
        Since $\varphi,\psi\in L^\infty(\mathbb{T}^2)$, let
        $$
        \varphi(z)=\sum_{m,n\in\mathbb{Z}}\widehat{\varphi}(m,n)z_1^mz_2^n,\quad\psi(z)=\sum_{k,l\in\mathbb{Z}}\widehat{\psi}(k,l)z_1^kz_2^l,
        $$
        where $\widehat{\varphi}(m,n)$ and $\widehat{\psi}(k,l)$ are the Fourier coefficients of $\varphi$ and $\psi$ in the $L^2(\mathbb T^2)$ sense respectively.
        Then for $(i,j)\in S$,
        $$
        \widehat{T}_\varphi\widehat{T}_\psi\chi_{i,j}=\sum_{(k+i,l+j)\in S,(m+k+i,n+l+j)\in S}\widehat{\varphi}(m,n)\widehat{\psi}(k,l)z_1^{m+k+i}z_2^{n+l+j}
        $$
        and
        $$
        \widehat{T}_\psi\widehat{T}_\varphi\chi_{i,j}=\sum_{(m+i,n+j)\in S,(m+k+i,n+l+j)\in S}\widehat{\varphi}(m,n)\widehat{\psi}(k,l)z_1^{m+k+i}z_2^{n+l+j}.
        $$
        Since 
        $$
        \widehat{T}_\varphi\widehat{T}_\psi\chi_{i,j}=\widehat{T}_\psi\widehat{T}_\varphi\chi_{i,j},\quad\forall (i,j)\in S,
        $$
        let $M=m+k+i$ and $N=n+l+j$, then we have for any $(M,N)\in S$,
        $$
        \sum_{(k+i,l+j)\in S}\widehat{\varphi}(M-i-k,N-j-l)\widehat{\psi}(k,l)=\sum_{(m+i,n+j)\in S}\widehat{\varphi}(m,n)\widehat{\psi}(M-i-m,N-j-n),
        $$
        i.e.
        $$
        (\widehat{\varphi}\ast_S^{(i,j)}\widehat{\psi})(M,N)=(\widehat{\psi}\ast_S^{(i,j)}\widehat{\varphi})(M,N),\quad\forall(i,j)\in S.
        $$
        Hence for any $(i,j)\in S$, we have
        $$
        \widehat{\varphi}\ast_S^{(i,j)}\widehat{\psi}=\widehat{\psi}\ast_S^{(i,j)}\widehat{\varphi}\quad\text{on $S$.}
        $$

        Conversely, assume that for any $(i,j)\in S$, we have
        $$
        \widehat{\varphi}\ast_S^{(i,j)}\widehat{\psi}=\widehat{\psi}\ast_S^{(i,j)}\widehat{\varphi}\quad\text{on $S$.}
        $$
        Fixing $(i,j)\in S$. By the Fourier expansions established above, we have
        $$
        \begin{aligned}
            \widehat{T}_\varphi\widehat{T}_\psi\chi_{i,j}&=\sum_{(M,N)\in S}\left[\sum_{(k+i,l+j)\in S}\widehat{\varphi}(M-i-k,N-j-l)\widehat{\psi}(k,l)\right]z_1^Mz_2^N\\
            &=\sum_{(M,N)\in S}(\widehat{\varphi}\ast_S^{(i,j)}\widehat{\psi})(M,N)z_1^Mz_2^N
        \end{aligned}
        $$
        and
        $$
        \begin{aligned}
            \widehat{T}_\psi\widehat{T}_\varphi\chi_{i,j}&=\sum_{(M,N)\in S}\left[\sum_{(m+i,n+j)\in S}\widehat{\varphi}(m,n)\widehat{\psi}(M-i-m,N-j-n)\right]z_1^Mz_2^N\\
            &=\sum_{(M,N)\in S}(\widehat{\psi}\ast_S^{(i,j)}\widehat{\varphi})(M,N)z_1^Mz_2^N.    
        \end{aligned}
        $$
        Hence
        $$
        \widehat{T}_\varphi\widehat{T}_\psi\chi_{i,j}=\widehat{T}_\psi\widehat{T}_\varphi\chi_{i,j},\quad\forall(i,j)\in S.
        $$
        Since $\{\chi_{i,j}\}_{(i,j)\in S}$
        is an orthonormal basis of $h^2(\mathbb{T}^2)$ and both operators are bounded, we have $\widehat{T}_\varphi$ commutes with $\widehat{T}_\psi$ on $h^2(\mathbb{T}^2)$.
        This completes the proof.
    \end{proof}

    For $\varphi\in L^2(\mathbb{T}^2)$, the Hankel operator $\widehat{H}_\varphi$ with symbol $\varphi$ is defined on the pluriharmonic Hardy
    space $h^2(\mathbb{T}^2)$ by
    $$
    \widehat{H}_\varphi:h^2(\mathbb{T}^2)\to\left(h^2(\mathbb{T}^2)\right)^\perp,\quad \widehat{H }_\varphi f=(I-Q)(\varphi f),\forall f\in h^2(\mathbb{T}^2).
    $$
    Define that
    $$
    \widehat{H}^*_\varphi:\left(h^2(\mathbb{T}^2)\right)^\perp\to h^2(\mathbb{T}^2),\quad\widehat{H}^*_\varphi f=Q(\bar{\varphi}f),\forall f\in\left(h^2(\mathbb{T}^2)\right)^\perp.
    $$

    For $f$ in $h^2(\mathbb{T}^2)$ and $\varphi$ in $L^2(\mathbb{T}^2)$, we have
    $$
    \varphi f=Q(\varphi f)+(I-Q)(\varphi f)=\widehat{T}_\varphi f+\widehat{H }_\varphi f.
    $$
    Then for all $\varphi,\psi$ in $L^\infty(\mathbb{T}^2)$ and $f$ in $h^\infty(\mathbb{T}^2)$,
    $$
    \varphi\psi f=\varphi(\widehat{T}_\psi f+\widehat{H}_\psi f)=\varphi\widehat{T }_\psi f+\varphi\widehat{H}_\psi f.
    $$
    It follows that
    $$
    \begin{aligned}
        \widehat{T}_{\varphi\psi}f&=Q(\varphi\psi f)\\
        &=Q(\varphi\widehat{T }_\psi f+\varphi\widehat{H}_\psi f)\\
        &=Q(\varphi\widehat{T}_\psi f)+Q(\varphi\widehat{H}_\psi f)\\
        &=\widehat{T}_\varphi\widehat{T}_\psi f+\widehat{H}^*_{\bar{\varphi}}\widehat{H}_\psi f,
    \end{aligned}
    $$
    since $\widehat{T}_\psi f$ is in $h^\infty(\mathbb{T}^2)$ and $\widehat{H}_\psi f$ is in $\left(h^\infty(\mathbb{T}^2)\right)^\perp$.
    Hence, we obtain that
    $$
    \widehat{T}_{\varphi\psi}=\widehat{T}_\varphi\widehat{T}_\psi+\widehat{H}^*_{\bar{\varphi}}\widehat{H }_\psi,
    $$
    for $\varphi,\psi$ in $L^\infty(\mathbb{T}^2)$.

    In the Fourier index plane, define the set
    $$
    S_{IV}=\{(u,v)\in\mathbb{Z}^2:u\geqslant0,v\leqslant0\}.
    $$

    In the product identity above, the Hankel term vanishes when both symbols have Fourier support in $S_{IV}$.
    Applying the identity in both orders gives the following corollary.

    \begin{coro}
        Let $\varphi,\psi\in L^\infty(\mathbb{T}^2)$. If $\varphi$ and $\psi$ satisfy
        $$
        \mathrm{supp}\widehat{\varphi}\subseteq S_{IV},\quad\mathrm{supp}\widehat{\psi}\subseteq S_{IV},
        $$
        then
        $$
        \widehat{T}_\varphi\widehat{T}_\psi=\widehat{T}_{\varphi\psi}=\widehat{T}_\psi\widehat{T}_\varphi.
        $$
    \end{coro}

    \begin{proof}[Proof]
        For any $f\in h^2(\mathbb{T}^2)$, its Fourier support lies in $S$. 
        Since $\mathrm{supp}\widehat{\psi}\subseteq S_{IV}$, $\psi f$ cannot have frequencies with
        $$
        u<0,\quad v>0,
        $$
        since if $i\geqslant0,j\geqslant0$, then $u=i+k\geqslant0$, and if $i\leqslant0,j\leqslant0$, then $v=j+l\leqslant0$.
        Therefore, the Fourier support of $(I-Q)(\psi f)$ must lie in $S_{IV}$.

        Then multiplying $\varphi$, since $\mathrm{supp}\widehat{\varphi}\subseteq S_{IV}$, these frequencies satisfy
        $$
        u+m>0,\quad v+n<0,
        $$
        and they are not in $S$. Hence
        $$
        Q[\varphi(I-Q)(\psi f)]=0,
        $$
        i.e. $\widehat{H}^*_{\bar{\varphi}}\widehat{H }_\psi f=0$.
        Therefore, by
        $$
        \widehat{T}_{\varphi\psi}=\widehat{T}_\varphi\widehat{T}_\psi+\widehat{H}^*_{\bar{\varphi}}\widehat{H }_\psi,
        $$
        we get
        $$
        \widehat{T}_\varphi\widehat{T}_\psi f=\widehat{T}_{\varphi\psi}f.
        $$
        Interchanging $\varphi$ and $\psi$ gives the same result.
        This completes this proof.
    \end{proof}

    \section{Commutativity of Toeplitz and little Hankel operators}

    In this section, we characterize the commutativity of Toeplitz and little Hankel operators on the pluriharmonic Hardy space of the bidisc. 
    We begin with the case of bounded symbols depending on a single variable similarly.

    \begin{thm}\label{thm3.13}
        Let $\varphi,\psi\in L^\infty(\mathbb{T})$. $\widehat{T}_\varphi$ commutes with $\widehat{\Gamma}_\psi$ on $h^2(\mathbb{T}^2)$ if and only if one of the following condition holds:

        \begin{enumerate}[label = \textup{(\alph*)}]
            \item $\psi=0$;
            \item $\varphi$ is a constant;
            \item $\widetilde{\varphi}=\varphi$, and there exists $c,d$ in $\mathbb{C}$ s.t. $\psi=c\varphi+d.$
        \end{enumerate}
    \end{thm}

    \begin{proof}[Proof]
        Since $S=\{(u,v)\in\mathbb{Z}^2:u\geqslant0,v\geqslant0\}\cup\{(u,v)\in\mathbb{Z}^2:u\leqslant0,v\leqslant0\}=-S$,
        and $U\chi_{u,v}=\chi_{-u,-v}$, we have $UQ=QU$ and $U^2=I$. It follows that
        $$
        \widehat{\Gamma}_\psi=U\widehat{T}_\psi,\quad U\widehat{T}_\varphi U=\widehat{T}_{\widetilde{\varphi}}.
        $$
        Hence we have
        $$
        \widehat{T}_\varphi\widehat{\Gamma}_\psi=\widehat{\Gamma}_\psi\widehat{T}_\varphi\quad\Leftrightarrow\quad\widehat{T}_\varphi U\widehat{T}_\psi=U\widehat{T}_\psi\widehat{T}_\varphi\quad\Leftrightarrow\quad\widehat{T}_{\widetilde{\varphi}}\widehat{T}_\psi=\widehat{T}_\psi\widehat{T}_\varphi.
        $$

        Since $\varphi$ and $\psi$ are one variable functions in $L^\infty{\mathbb{T}}$, by commutativity, we can clearly get $\psi=0$ or $\varphi$ is a constant.
        Now we assume that $\psi\neq0$ and $\varphi$ is non-constant. 
        
        Firstly we assume that $\varphi$ and $\psi$ have the same variable $z_1$, i.e.
        $$
        \varphi(z)=\sum_{m\in\mathbb{Z}}\widehat{\varphi}(m)z_1^m,\quad\psi(z)=\sum_{n\in\mathbb{Z}}\widehat{\psi}(n)z_1^n,
        $$
        where $\widehat{\varphi}(m)$ and $\widehat{\psi}(n)$ are the Fourier coefficients of $\varphi$ and $\psi$ in the $L^2(\mathbb T^2)$ sense respectively.
        Since $\widehat{T}_\varphi\widehat{\Gamma}_\psi\chi_{k,j}=\widehat{\Gamma}_\psi\widehat{T}_\varphi\chi_{k.j}$,
        for all $(k,j)\in S$, by calculating, we can get
        $$
        \widehat{T}_\varphi\widehat{\Gamma}_\psi\chi_{k,j}=\begin{cases}
            \displaystyle\sum_{n+k\geqslant0,m-n-k\leqslant0}\widehat{\varphi}(m)\widehat{\psi}(n)z_1^{m-n-k}z_2^{-j},\quad&k\geqslant0,j>0,\\
            \vspace{2mm}\\
            \displaystyle\sum_{n+k\leqslant0,m-n-k\geqslant0}\widehat{\varphi}(m)\widehat{\psi}(n)z_1^{m-n-k}z_2^{-j},&k\leqslant0,j<0,\\
            \vspace{2mm}\\
            \displaystyle\sum_{m,n\in\mathbb{Z}}\widehat{\varphi}(m)\widehat{\psi}(n)z_1^{m-n-k},&k\in\mathbb{Z},j=0
        \end{cases}
        $$
        and
        $$
        \widehat{\Gamma}_\psi\widehat{T}_\varphi\chi_{k,j}=\begin{cases}
            \displaystyle\sum_{m+k\geqslant0,m+n+k\geqslant0}\widehat{\varphi}(m)\widehat{\psi}(n)z_1^{-(m+n+k)}z_2^{-j},\quad&k\geqslant0,j>0,\\
            \vspace{2mm}\\
            \displaystyle\sum_{m+k\leqslant0,m+n+k\leqslant0}\widehat{\varphi}(m)\widehat{\psi}(n)z_1^{-(m+n+k)}z_2^{-j},&k\leqslant0,j<0,\\
            \vspace{2mm}\\
            \displaystyle\sum_{m,n\in\mathbb{Z}}\widehat{\varphi}(m)\widehat{\psi}(n)z_1^{-(m+n+k)},&k\in\mathbb{Z},j=0.
        \end{cases}
        $$

        For $k\in\mathbb{Z},j=0$, $Q(\chi_{k,0})=\chi_{k,0}$. So taking $k=0,j=0$ we have
        $$
        \varphi\widetilde{\psi}=\widetilde{\varphi}\widetilde{\psi}.
        $$
        Applying $U$ to both sides of the equation yields
        \begin{equation}\label{3.12.1}
            \psi(\varphi-\widetilde{\varphi})=0\quad\text{a.e.}
        \end{equation}

        Then for $k\geqslant0,j>0$, we have
        $$
        \sum_{n+k\geqslant0,m-n-k\leqslant0}\widehat{\varphi}(m)\widehat{\psi}(n)z_1^{m-n-k}z_2^{-j}=\sum_{m+k\geqslant0,m+n+k\geqslant0}\widehat{\varphi}(m)\widehat{\psi}(n)z_1^{-(m+n+k)}z_2^{-j}.
        $$
        Let $A=m-n-k$ and $B=-m-n-k$,
        $$
        \sum_{n+k\geqslant0,A\leqslant0}\widehat{\varphi}(m)\widehat{\psi}(n)z_1^Az_2^{-j}=\sum_{m+k\geqslant0,B\leqslant0}\widehat{\varphi}(m)\widehat{\psi}(n)z_1^Bz_2^{-j}.
        $$
        For $A=B\leqslant0$, we can get
        $$
        \sum_{n+k\geqslant0}\widehat{\varphi}(A+n+k)\widehat{\psi}(n)=\sum_{m+k\geqslant0}\widehat{\varphi}(m)\widehat{\psi}(-B-m-k).
        $$
        Let $-q=A=B\leqslant0$ and $s=n+k$ in left-hand side, $s=m+k$ in right-hand side, we have
        \begin{equation}\label{3.12.2}
            \sum_{s\geqslant0}\widehat{\varphi}(s-q)\widehat{\psi}(s-k)=\sum_{s\geqslant0}\widehat{\varphi}(s-k)\widehat{\psi}(q-s).
        \end{equation}
        Similarly as the case $k\leqslant0,j<0$, we can get
        \begin{equation}\label{3.12.3}
            \sum_{s\geqslant0}\widehat{\varphi}(q-s)\widehat{\psi}(-s-k)=\sum_{s\geqslant0}\widehat{\varphi}(-s-k)\widehat{\psi}(s-q),\quad q\geqslant0,k\leqslant0.
        \end{equation}
        Here, the corresponding series is absolutely convergent for every fixed $k$ and $q$.
        Replacing $q$ and $k$ in (\ref{3.12.2}) simultaneously by $q+1$ and $k+1$, then changing the
        summation index to $s-1$, we have
        $$
        \widehat{\varphi}(-q-1)\widehat{\psi}(-k-1)+\sum_{s\geqslant0}\widehat{\varphi}(s-q)\widehat{\psi}(s-k)=\widehat{\varphi}(-k-1)\widehat{\psi}(q+1)+\sum_{s\geqslant0}\widehat{\varphi}(s-k)\widehat{\psi}(q-s),
        $$
        and subtracting the result from the original equation, we obtain
        $$
        \widehat{\varphi}(-q-1)\widehat{\psi}(-k-1)=\widehat{\varphi}(-k-1)\widehat{\psi}(q+1).
        $$
        Let $m=q+1,n=k+1$, then
        \begin{equation}\label{3.12.4}
            \widehat{\varphi}(-m)\widehat{\psi}(-n)=\widehat{\varphi}(-n)\widehat{\psi}(m),\quad m\geqslant1,n\geqslant1.
        \end{equation}
        Allpying the same method to (\ref{3.12.3}) yields
        \begin{equation}\label{3.12.5}
            \widehat{\varphi}(m)\widehat{\psi}(n)=\widehat{\varphi}(n)\widehat{\psi}(-m),\quad m\geqslant1,n\geqslant1.
        \end{equation}

        Since $\varphi$ is non-constant, there exists $p\geqslant1$ s.t. $\widehat{\varphi}(p)\neq0$ or $\widehat{\varphi}(-p)\neq0$.
        If $\widehat{\varphi}(p)\neq0$, let $m=n=p$ in (\ref{3.12.5}), then
        \begin{equation}\label{3.12.6}
            \widehat{\psi}(p)=\widehat{\psi}(-p).
        \end{equation}
        Then, setting $(m,n)=(p,r)$ and $(r,p)$, $r\geqslant1$, in (\ref{3.12.5}), respectively, we obtain
        $$
        \widehat{\varphi}(p)\widehat{\psi}(r)=\widehat{\varphi}(r)\widehat{\psi}(-p),\quad\widehat{\varphi}(r)\widehat{\psi}(p)=\widehat{\varphi}(p)\widehat{\psi}(-r).
        $$
        Hence, by (\ref{3.12.6}), we have
        $$
        \widehat{\psi}(r)=\widehat{\psi}(-r),\quad r\geqslant1.
        $$
        If $\widehat{\varphi}(-p)\neq0$, using (\ref{3.12.4}) in the same way also yields this conclusion.
        Hence, we have
        \begin{equation}\label{3.12.7}
            \widetilde{\psi}=\psi.
        \end{equation}

        If $\psi$ is a constant, since $\psi\neq0$ and (\ref{3.12.1}), we have
        $$
        \widetilde{\varphi}=\varphi.
        $$
        Then, with $c=0$ and $d=\psi$, condition (c) holds. Next consider $\psi$ is non-constant.
        By (\ref{3.12.7}), there exists $p\geqslant1$ s.t. $\widehat{\psi}(p)=\widehat{\psi}(-p)\neq0$.
        Define two constants
        $$
        c_1=\frac{\widehat{\varphi}(p)}{\widehat{\psi}(p)},\quad c_2=\frac{\widehat{\varphi}(-p)}{\widehat{\psi}(p)}.
        $$
        By (\ref{3.12.4}) and (\ref{3.12.5}), for any $r\geqslant1$,
        \begin{equation}\label{3.12.8}
            \widehat{\varphi}(r)=c_1\widehat{\psi}(r),\quad\widehat{\varphi}(-r)=c_2\widehat{\psi}(r).
        \end{equation}
        Now we prove that $c_1=c_2$. Assume that $c_1\neq c_2$, and define function
        $$
        f=\frac{\varphi-\widehat{\varphi}(0)-c_2(\psi-\widehat{\psi}(0))}{c_1-c_2}.
        $$
        Obviously $f$ is in $L^\infty(\mathbb{T})$. By (\ref{3.12.7}) and (\ref{3.12.8}), we have
        $$
        \widehat{f}(0)=0,\quad\widehat{f}(r)=\widehat{\psi}(r),\quad\widehat{f}(-r)=0,\quad r\geqslant1.
        $$
        Hence $f$ is in $H^\infty(\mathbb{T})$ and $f(0)=0$, and
        $$
        \psi=\widehat{\psi}(0)+f+\widetilde{f},\quad\varphi=\widehat{\varphi}(0)+c_1f+c_2\widetilde{f}.
        $$
        Substituting there relations into (\ref{3.12.1}) and conceling the nonzero constant $c_1-c_2$, we obtain
        $$
        (f-\widetilde{f})(\widehat{\psi}(0)+f+\widetilde{f})=0,
        $$
        i.e.
        $$
        f^2+\widehat{\psi}(0)f=\widetilde{f}^2+\widehat{\psi}(0)\widetilde{f}=\widetilde{f^2+\widehat{\psi}(0)f}.
        $$
        The left-hand side lies in $H^\infty(\mathbb{T})$, while the right-hand side has lies in $\overline{H^\infty(\mathbb{T})}$.
        Hence the equality of the two sides forces the function to be constant. Since $f(0)=0$, this constants is zero, i.e.
        $$
        f(f+\widehat{\psi}(0))\equiv0.
        $$
        Since $f$ is analytic, we have $f\equiv0$. But $\widehat{f}(p)=\widehat{\psi}(p)\neq0$. Contradiction!
        Hence $c_1=c_2$. Then by (\ref{3.12.8}), we have
        $$
        \widetilde{\varphi}=\varphi,\quad\varphi-\widetilde{\varphi}(0)=c_1(\psi-\widehat{\psi}(0)).
        $$
        Since $\varphi$ is non-constant, $c_1\neq0$. Taking
        $$
        c=\frac{1}{c_1},\quad d=\widehat{\psi}(0)-\frac{\widehat{\varphi}(0)}{c_1},
        $$
        we have
        $$
        \psi=c\varphi+d.
        $$
            
        For $\varphi$ and $\psi$ depending on $z_2$, the argument proceeds in exactly the same way as above. Therefore (c) holds.

        If $\varphi$ and $\psi$ are in different variable, we may assume that $\varphi$ depends on $z_1$ and $\psi$ depends on $z_2$, i.e.
        $$
        \varphi(z)=\sum_{m\in\mathbb{Z}}\widehat{\varphi}(m)z_1^m,\quad\psi(z)=\sum_{n\in\mathbb{Z}}\widehat{\psi}(n)z_2^n,
        $$
        and $\psi\neq0$ and $\varphi$ is non-constant. For $k>0,j>0$,
        $$
        \widehat{T}_{\widetilde{\varphi}}\widehat{T}_\psi\chi_{k,j}=\sum_{n+j\geqslant0,k-m\geqslant0}\widehat{\varphi}(m)\widehat{\psi}(n)z_1^{k-m}z_2^{n+j}+\widehat{\psi}(-j)\sum_{k-m<0}\widehat{\varphi}(m)z_1^{k-m}
        $$
        and
        $$
        \widehat{T}_\psi\widehat{T}_\varphi\chi_{k,j}=\sum_{m+k\geqslant0,n+j\geqslant0}\widehat{\varphi}(m)\widehat{\psi}(n)z_1^{m+k}z_2^{n+j}+\widehat{\varphi}(-k)\sum_{n+j<0}\widehat{\psi}(n)z_2^{n+j}.
        $$
        Since $\widehat{T}_{\widetilde{\varphi}}\widehat{T}_\psi=\widehat{T}_\psi\widehat{T}_\varphi$, we have
        \begin{equation}\label{3.12.9}
            \sum_{n+j\geqslant0,k-m\geqslant0}\widehat{\varphi}(m)\widehat{\psi}(n)z_1^{k-m}z_2^{n+j}=\sum_{m+k\geqslant0,n+j\geqslant0}\widehat{\varphi}(m)\widehat{\psi}(n)z_1^{m+k}z_2^{n+j}
        \end{equation}
        and
        \begin{equation}\label{3.12.10}
            \widehat{\psi}(-j)\sum_{k-m<0}\widehat{\varphi}(m)z_1^{k-m}=\widehat{\varphi}(-k)\sum_{n+j<0}\widehat{\psi}(n)z_2^{n+j}.
        \end{equation}
        Obviously, (\ref{3.12.10}) must be zero. Then for (\ref{3.12.9}), let $q=k-m,s=n+j$ in left-hand side and $q=m+k,s=n+j$ in right-hand side, then
        $$
        \sum_{q,s\geqslant0}\widehat{\varphi}(k-q)\widehat{\psi}(s-j)z_1^qz_2^s=\sum_{q,s\geqslant0}\widehat{\varphi}(q-k)\widehat{\psi}(s-j)z_1^qz_2^s.
        $$
        It follows that
        $$
        [\widehat{\varphi}(k-q)-\widehat{\varphi}(q-k)]\widehat{\psi}(s-j)=0,\quad k>0,j>0,q\geqslant,s\geqslant0.
        $$
        Since $\psi\neq0$, there exists $n_0\in\mathbb{Z}$ s.t. $\widehat{\psi}(n_0)\neq0$.
        Let $j>\max\{0,-n_0\}$ and $s=j+n_0>0$, then
        $$
        \widehat{\varphi}(k-q)=\widehat{\varphi}(q-k),\quad k>0,q\geqslant0,
        $$
        where $q-k$ is in $\mathbb{Z}$. It means that for any $m$ in $\mathbb{Z}$, let $k>\max\{0,-m\}$ and $q=k+m>0$, then
        $$
        \widehat{\varphi}(-m)=\widehat{\varphi}(m),
        $$
        i.e. $\widetilde{\varphi}=\varphi$. Thus, the original operator equation becomes
        $$
        \widehat{T}_\varphi\widehat{T}_\psi=\widehat{T}_\psi\widehat{T}_\varphi.
        $$
        By Theorem \ref{thm3.10}, if $\varphi$ and $\psi$ are in different variable, then one of $\varphi$ and $\psi$ is analytic and the other is co-analytic.
        But we have known $\widetilde{\varphi}=\varphi$, i.e. if $\varphi$ is analytic or co-analytic, $\varphi$ must be a constant. Contradiction!
        Hence $\psi$ must be a constant by Theorem \ref{thm3.10}. Setting $c=0$ and $d=\psi$, we obtain $\psi=c\varphi+d$, which completes the proof of condition (c).

        Conversely, if condition (a) or (b) holds, the commutativity is clearly.
        For condition (c), since $\widetilde{\varphi}=\varphi$, we have
        $$
        U\widehat{T}_\varphi f=UQ(\varphi f)=QU(\varphi f)=Q(\widetilde{\varphi}Uf)=Q(\varphi Uf)=\widehat{T}_\varphi U f,\quad\forall f\in h^2(\mathbb{T}^2).
        $$
        Since $\psi=c\varphi+d$,
        $$
        \widehat{\Gamma}_\psi=U(c\widehat{T}_\varphi+dI).
        $$
        Hence
        $$
        \widehat{T}_\varphi\widehat{\Gamma}_\psi=\widehat{T}_\varphi U(c\widehat{T}_\varphi+dI)=U(c\widehat{T}_\varphi+dI)\widehat{T}_\varphi=\widehat{\Gamma}_\psi\widehat{T}_\varphi.
        $$
        This completes the proof.
    \end{proof}

    We next consider the same commutation problem for general symbols in $L^\infty(\mathbb{T}^2)$.

    \begin{thm}\label{thm3.14}
        Let $\varphi,\psi\in L^\infty(\mathbb{T}^2)$. Then $\widehat{T}_\varphi$ commutes with $\widehat{\Gamma}_\psi$ on $h^2(\mathbb{T}^2)$ if and only
        if for every $(i,j)\in S$,
        $$
        \widehat{\widetilde{\varphi}}\ast_S^{(i,j)}\widehat{\psi}=\widehat{\psi}\ast_S^{(i,j)}\widehat{\varphi}\quad\text{on $S$.}
        $$
    \end{thm}

    \begin{proof}[Proof]
        By the proof of Theorem \ref{thm3.13}, we have known
        \begin{equation}\label{3.14.1}
            \widehat{T}_\varphi\widehat{\Gamma}_\psi=\widehat{\Gamma}_\psi\widehat{T}_\varphi\quad\Leftrightarrow\quad\widehat{T}_{\widetilde{\varphi}}\widehat{T}_\psi=\widehat{T}_\psi\widehat{T}_\varphi.
        \end{equation}
        Fix $(i,j)\in S$. By the product coefficient formulas established in the proof of Theorem \ref{thm3.11},,
        $$
        \begin{aligned}
            \widehat{T}_{\widetilde{\varphi}}\widehat{T}_\psi\chi_{i,j}&=\sum_{(M,N)\in S}\left[\sum_{(k+i,l+j)\in S}\widehat{\widetilde{\varphi}}(M-i-k,N-j-l)\widehat{\psi}(k,l)\right]z_1^Mz_2^N\\
            &=\sum_{(M,N)\in S}(\widehat{\widetilde{\varphi}}\ast_S^{(i,j)}\widehat{\psi})(M,N)z_1^Mz_2^N
        \end{aligned}
        $$
        and
        $$
        \begin{aligned}
            \widehat{T}_\psi\widehat{T}_\varphi\chi_{i,j}&=\sum_{(M,N)\in S}\left[\sum_{(m+i,n+j)\in S}\widehat{\varphi}(m,n)\widehat{\psi}(M-i-m,N-j-n)\right]z_1^Mz_2^N\\
            &=\sum_{(M,N)\in S}(\widehat{\psi}\ast_S^{(i,j)}\widehat{\varphi})(M,N)z_1^Mz_2^N.    
        \end{aligned}
        $$
        
        Suppose first that $\widehat{T}_\varphi$ commutes with $\widehat{\Gamma}_\psi$. By (\ref{3.14.1}),
        $$
        \widehat{T}_{\widetilde{\varphi}}\widehat{T}_\psi\chi_{i,j}=\widehat{T}_\psi\widehat{T}_\varphi\chi_{i,j}
        $$
        for every $(i,j)\in S$. Comparing the coefficients gives
        $$
        (\widehat{\widetilde{\varphi}}\ast_S^{(i,j)}\widehat{\psi})(M,N)=(\widehat{\psi}\ast_S^{(i,j)}\widehat{\varphi})(M,N)
        $$
        for every $(M,N)\in S$. Therefore, we have for every $(i,j)\in S$,
        $$
        \widehat{\widetilde{\varphi}}\ast_S^{(i,j)}\widehat{\psi}=\widehat{\psi}\ast_S^{(i,j)}\widehat{\varphi}\quad\text{on $S$.}
        $$

        Conversely, assume that
        $$
        \widehat{\widetilde{\varphi}}\ast_S^{(i,j)}\widehat{\psi}=\widehat{\psi}\ast_S^{(i,j)}\widehat{\varphi}\quad\text{on $S$}
        $$
        for every $(i,j)\in S$. The preceding Fourier expansions yields that
        $$
        \widehat{T}_{\widetilde{\varphi}}\widehat{T}_\psi\chi_{i,j}=\widehat{T}_\psi\widehat{T}_\varphi\chi_{i,j},\quad\forall(i,j)\in S.
        $$
        Since $\{\chi_{i,j}\}_{(i,j)\in S}$
        is an orthonormal basis of $h^2(\mathbb{T}^2)$ and both operators are bounded, 
        we have $\widehat{T}_{\widetilde{\varphi}}\widehat{T}_\psi=\widehat{T}_\psi\widehat{T}_\varphi$ on $h^2(\mathbb{T}^2)$,
        i.e. $\widehat{T}_\varphi$ commutes with $\widehat{\Gamma}_\psi$ on $h^2(\mathbb{T}^2)$.
        This completes the proof.
    \end{proof}

    \subsection*{Acknowledgment}
    This work was supported by the National Natural Science Foundation of China (12371125) and the Chongqing Natural Science Foundation (CSTB2024NSCQ-MSX0177).


\begin{thebibliography}{99}

        \bibitem{BH1964}
        A. Brown, P. R. Halmos,
        Algebraic properties of Toeplitz operators[J],
        Journal f\"{u}r die reine und angewandte Mathematik,
        1964, 213: 89--102.

        \bibitem{AC1991}
        S. Axler, \v{Z}. \v{C}u\v{c}kovi\'{c},
        Commuting Toeplitz operators with harmonic symbols[J],
        Integral Equations and Operator Theory,
        1991, 14(1): 1--12.

        \bibitem{AC2001}
        P. Ahern, \v{Z}. \v{C}u\v{c}kovi\'{c},
        A theorem of Brown--Halmos type for Bergman space Toeplitz operators[J],
        Journal of Functional Analysis,
        2001, 187(1): 200--210.

        \bibitem{GZ1997}
        C. Gu, D. Zheng,
        The semi-commutator of Toeplitz operators on the bidisc[J],
        Journal of Operator Theory,
        1997, 38(1): 173--193.

        \bibitem{Ding2003}
        X. Ding,
        Products of Toeplitz operators on the polydisk[J],
        Integral Equations and Operator Theory,
        2003, 45(4): 389--403.

        \bibitem{CKL2011}
        B. R. Choe, H. Koo, Y. J. Lee,
        Toeplitz products with pluriharmonic symbols on the Hardy space over the ball[J],
        Journal of Mathematical Analysis and Applications,
        2011, 381(1): 365--382.

        \bibitem{DSZ2012}
        X. Ding, S. Sun, D. Zheng,
        Commuting Toeplitz operators on the bidisk[J],
        Journal of Functional Analysis,
        2012, 263(11): 3333--3357.

        \bibitem{KL2015}
        L. H. Kong, Y. F. Lu,
        Commuting Toeplitz operators on the Hardy space of the polydisk[J],
        Acta Mathematica Sinica, English Series,
        2015, 31(4): 695--702.

        \bibitem{DWL2023}
        X. Ding, Z. Wang, Y. Li,
        Some questions about the Hardy-Toeplitz operator of the bidisk[J],
        Scientia Sinica Mathematica,
        2023, 53(10): 1349--1356.

        \bibitem{LD2013}
        Y. Liu, X. Ding,
        Toeplitz operators on the pluriharmonic Hardy space[J],
        Scientia Sinica Mathematica,
        2013, 43(7): 675--690.

        \bibitem{DS2017}
        X. Ding, Y. Sang,
        The pluriharmonic Hardy space and Toeplitz operators[J],
        Results in Mathematics,
        2017, 72(3): 1473--1497.

        \bibitem{SD2018}
        Y. Sang, X. Ding,
        Commuting Toeplitz operators on the pluriharmonic Hardy space[J],
        Acta Mathematica Sinica, Chinese Series,
        2018, 61(4): 577--584.

        \bibitem{ZD2020}
        Y. Zhang, X. Dong,
        Commutators and semi-commutators of monomial Toeplitz operators
        on the pluriharmonic Hardy space[J],
        Chinese Annals of Mathematics, Series B,
        2020, 41(5): 717--732.

        \bibitem{MA2000}
        R. A. Mart\'inez-Avenda\~{n}o,
        When do Toeplitz and Hankel operators commute?[J],
        Integral Equations and Operator Theory,
        2000, 37: 341--349.

        \bibitem{ZY2026}
        H. Zou, T. Yu,
        Commuting Hankel and Toeplitz operators on the Hardy space of the bidisk[J],
        Advances in Operator Theory,
        2026, 11: 41.
	\end{thebibliography}
\end{document}